\documentclass[12pt,reqno]{amsart}
\usepackage[utf8]{inputenc}
\usepackage{amsmath, amsfonts, amssymb, amsthm, hyperref}
\usepackage{bm}
\allowdisplaybreaks[4]
\def\l{\left}
\def\r{\right}
\def\bg{\bigg}
\def\({\bg(}
\def\){\bg)}
\def\t{\text}

\def\p{\mathfrak p}
\def\1{{\bf 1}}

\def\pmod #1{\ ({\rm{mod}}\ #1)}

\def\<{\langle}
\def\>{\rangle}

\theoremstyle{plain}
\newtheorem{theorem}{Theorem}[section]
\newtheorem{lemma}{Lemma}

\theoremstyle{definition}

\theoremstyle{remark}
\newtheorem{remark}{Remark}

\makeatletter
\@namedef{subjclassname@2024}{%
  \textup{2024} Mathematics Subject Classification}
\makeatother
 \vspace{4mm}

\begin{document}
\medskip

\title[On the natural density of
 integers $n$  for which $\sigma(kn+r_1) >\sigma(kn+r_2)$ ]
{ On the natural density
of integers $n$ \\ for which $\sigma(kn+r_1) >\sigma(kn+r_2)$ }
\author{Xin-qi Luo}
\address {(Xin-qi Luo) School of Sciences, Changzhou institute of technology, Changzhou 213032, People's Republic of China}
\email{luoxq@czust.edu.cn}

\author{Chen-kai Ren}
\address {(Chen-kai Ren) Department of Mathematics, Nanjing
University, Nanjing 210093, People's Republic of China}
\email{ckren@smail.nju.edu.cn}

\keywords{the sum of divisors function, density estimate, partition.
\newline \indent Supported by the Natural Science Foundation of China (grant no. 12371004).}
\begin{abstract}
For any positive integer $n$, let $\sigma(n)=\sum_{d\mid n} d$. In 2020, M. Kobayashi and T. Trudgian showed that the natural density of positive integers n with $\sigma(kn+r_1) \geq \sigma(kn+r_2)$ is between 0.053 and 0.055. In this paper, we extend their result. For integers $k>r_1>r_2\geq 0,$ we provide an estimate on the natural density of positive integers $n$ for which $\sigma(kn+r_1) > \sigma(kn+r_2)$. We also compute explicit bounds for specific $k,r_1,r_2$ to illustrate the variation of density.

\end{abstract}
\maketitle

\section{Introduction}
\setcounter{lemma}{0}
\setcounter{theorem}{0}
\setcounter{equation}{0}
\setcounter{conjecture}{0}
\setcounter{remark}{0}
\setcounter{corollary}{0}
For any positive integer $n$, let $\sigma(n) $ be the sum of all positive divisors of $n$. There are many interesting results on this arithmetical function. In 1936, Erd\"os \cite{E36} proved that the natural density of the set of positive integers $n$ for which $\sigma(n+1) \geq \sigma(n)$ is $\frac{1}{2}$. And in 1980, Erd\"os, Gy\"ory and Papp \cite{E80} proved that the chain of inequalities
$$
\sigma(n+m_1)>\sigma(n+m_2)>\sigma(n+m_3)>\sigma(n+m_4)
$$
holds for infinitely many $n$, where the $m_i$ are any permutations of the numbers $1,2,3,4$.
In 2022, R.-J. Wang and Y.-G. Chen \cite{WC22} showed that for  any integer $l\geq 2$,  $\sigma_{l}(2n+1)<\sigma_{l}(2n)$ and $\sigma_{l}(2n-1)<\sigma_{l}(2n)$ hold for all sufficiently large integers $n$, where $ \sigma_{l}(n)=\sum_{d\mid n}{d^l}.$

 Let $B$ be the set of natural numbers $n$ satisfying $\sigma(2n+1)\geq \sigma(2n)$ and $B(x)$ be the number of those $n$ in $B$ with $n\leq x$.
In 1987, Laub \cite{L87} posed the question of estimating the size of $B(x)/x$. Mattics \cite{M90} showed that there exist constants $\lambda$ and $\mu$ with
 $0 < \lambda < \mu < 1$ such that
 $$\lambda x < B(x) < \mu x$$
 for all large integers $x$ and recorded a remark of Hildebrand that
$\lim_{x\to \infty} B(x)/x$ exists. We call this limit the natural density of $B$, denoted $\mathbf{d} B.$    In 2020, M. Kobayashi and T. Trudgain \cite{KT20} refined Mattics' result and obtained  $0.0539171\leq \mathbf{d} B\leq 0.0549445$.

In this paper, we generalize Kobayashi and Trudgain's result and calculate some specific examples by computer.  Let $k,r_1$ and $r_2$ be  integers with $k>r_1>r_2\geq 0$ and let $A(k,r_1,r_2)$ be the set of natural numbers $n$ satisfying $\sigma(kn+r_1)>\sigma(kn+r_2)$. Then we have the following theorems:
 \begin{theorem}
 The density of $A(k,r_1,r_2)$ exists, and denoted $\mathbf{d}A(k,r_1,r_2)$.
 \end{theorem}

\begin{theorem}For $(k,r_1,r_2)=(3,2,0),$ we have
$$ 0.0361784 \leq \mathbf{d}A(3,2,0) \leq 0.1784994.  $$
 \end{theorem}
\begin{theorem}For $(k,r_1,r_2)=(4,1,0),$ we have
$$ 0.0107553 \leq \mathbf{d}A(4,1,0) \leq 0.0207187.  $$
 \end{theorem}
\begin{remark}\label{remark1}
  P. Erd\"os  \cite{E36} said that the density of integers $n$ with $\sigma (n)=\sigma (n+1)$ is zero. However, we met some difficulties in determining the density of integers $n$ with $\sigma(kn+r_1) = \sigma(kn+r_2).$ 
  We used the method developed in \cite{E36,P77} and we were left with the following case that we cannot handle:

  Let $l=\exp\{(\log x)^{\frac 18}\}$, $L=\exp\{\frac 18(\log x)^{\frac 23}\log \log x\}$ and $P(n)$ denote the largest prime factor of $n.$ For $i=1,2$, let $m_i=\frac {(kn+r_i)}{P(kn+r_i)}$. Then the number of $n<x$ satisfying the following conditions is $o(x)$:

  \noindent(1) $\sigma(kn+r_1)=\sigma(kn+r_2).$

  \noindent(2) If $t^a$ divides $kn+r_1$ or $kn+r_2$ with $a\in \{2,3,\ldots \}$, then $t^a\leq l^3$.

  \noindent(3) $P(kn+r_1)>L^2$ and $P(kn+r_2)>L^2.$

  \noindent(4) $m_1<L$ or $m_2<L.$

  \noindent(5)
  $m_1-m_2=r_1-r_2$ and $\sigma(m_1)/m_1=\sigma(m_2)/m_2.$

 \noindent Even with the help of A. Hildebrand's work \cite{H86}, we still cannot overcome this problem.
\end{remark}
\begin{remark}\label{remark2}
The upper and lower bounds are affected by the irregular summation terms which we will discuss in the third section, whose complexity introduces uncertainty.
\end{remark}
We are going to prove Theorem 1.1 in section 2. In section 3, we will provide an estimate of $\mathbf{d}A(k,r_1,r_2).$ And we complete the proof of Theorem 1.2 and 1.3 in section 4.

\section{Proof of Theorem 1.1 }
\setcounter{lemma}{0}
\setcounter{theorem}{0}
\setcounter{equation}{0}
\setcounter{conjecture}{0}
\setcounter{remark}{0}
\setcounter{corollary}{0}
\proof

Let $h(n)=\sigma(n)/n$. We define  $B(k,r_1,r_2)$ as
$$
B(k,r_1,r_2):=\{n:h(kn+r_1)\geq h(kn+r_2)\}.
$$
 By the equation

$$
\frac {h(kn+r_1)}{h(kn+r_2)}=\frac{\sigma (kn+r_1)}{\sigma (kn+r_2)}\cdot \frac {kn+r_2}{kn+r_1},
$$
we can find that $B(k,r_1,r_2)\subset A(k,r_1,r_2)$. Observe that
$$
A(k,r_1,r_2)-B(k,r_1,r_2)=\left\{n:0< \sigma(kn+r_1)-\sigma (kn+r_2)<\frac{(r_1-r_2)\sigma (kn+r_2)}{kn+r_2} \right\}.
$$
We will show that the set $A(k,r_1,r_2)-B(k,r_1,r_2)$
has density zero. By Gr\'onwall's theorem \cite{G13}, we obtain
$$ \limsup_{n\to \infty} \frac{\sigma(n)/n}{\log{\log{n}}}=e^\gamma,$$
where $\gamma$ is the Euler-Mascheroni constant. Therefore, for $n\in A(k,r_1,r_2)-B(k,r_1,r_2)$ we have
$$ \sigma(kn+r_1)-\sigma(kn+r_2)= O(\log{\log{n}}). $$
By Lemma 2.1 of \cite{LP15}, on a set $S$ of asymptotic density 1, $p \mid \sigma(n)$ for every prime $p\leq \log{\log{n}}/\log{\log{\log{n}}} $. Let $F(n)$ be the product of the primes not exceeding $\log{\log{n}}/\log{\log{\log{n}}}$. Since $\log{\log{(kn+r_2)}}/\log{\log{\log{(kn+r_2)}}}< \log{\log{(kn+r_1)}}/\log{\log{\log{(kn+r_1)}}}$, we have
$$ F(kn+r_2)\mid  \sigma(kn+r_1) \quad \it{and}  \quad F(kn+r_2)\mid  \sigma(kn+r_2).$$
Thus, for almost all $n$, we obtain that $ F(kn+r_2)\mid  \sigma(kn+r_1) - \sigma(kn+r_2). $ By the prime number theorem, we have
$$ F(n)=\log{(n)}^{(1+o(1))/\log{\log{\log{n}}}}. $$
 Hence, in set $A(k,r_1,r_2)-B(k,r_1,r_2)$, we have $\sigma(kn+r_1) >\sigma(kn+r_2)$ and
$$ \log{(kn+r_2)}^{(1+o(1))/\log{\log{\log{(kn+r_2)}}}} =F(kn+r_2)\leq \sigma(kn+r_1)-\sigma(kn+r_2)=  O(\log{\log{n}}), $$
which is a contradiction for sufficiently large $n$.

Since $B(k,r_1,r_2)$ has a density by \cite{S73}, we obtain that the set $A(k,r_1,r_2)$ has a density denoted $\mathbf{d}A(k,r_1,r_2)$ and the set $A(k,r_1,r_2)-B(k,r_1,r_2)$ has density zero. Now we complete the proof of Theorem 1.1. \qed

\section{ The estimate of bounds by partition }
\setcounter{lemma}{0}
\setcounter{theorem}{0}
\setcounter{equation}{0}
\setcounter{conjecture}{0}
\setcounter{remark}{0}
\setcounter{corollary}{0}

Let $y\geq 2.$ We say a number $n$ is $y-$smooth if its largest prime divisor $p$ satisfies $p\leq y$. Write $S(y)$ for the set of $y-$smooth numbers. Let $Y_y(n)$ be the largest $y-$smooth divisor of $n$. For integers $k>r_1>r_2\geq 0$, we define
$$
S_y(a,b):=\{n\in \mathbb N:Y_y(kn+r_1)=a, Y_y(kn+r_2)=b\}.
$$
We find that the sets $S_y(a,b)$ with $a,b\in S(y)$ partition $\mathbb{N}$ and that $S_y(a,b)=\emptyset $ unless  $\gcd(a,b)\mid r_1-r_2$.
We partition $B(k, r_1, r_2)$ via $B_y(a,b):=B(k, r_1, r_2)\cap S_y(a,b)$.
We denote the set of totatives modulo $N$ by
$$
\Phi(N):=\{t\in \mathbb N:1\leq t\leq N,\gcd(t,N)=1 \}.
$$
We define $P(y)$ as the product of all primes $p$ with $p\leq y.$ For any $n\in N,$ since $Y_y(n)$ is the largest $y-$smooth divisor of $n$ ,  $n/Y_y(n)$ does not have any prime divisor of $n$ not exceeding $y$. It follows that  $(n/Y_y(n),P(y))=1$. Hence, we may
partition $S_y(a,b)$ by
$$
S_y(a,b,t_1,t_2):=\{n\in S_y(a,b):(kn+r_1)/a\equiv t_1\pmod {P(y)},\ \ (kn+r_2)/b\equiv t_2 \pmod {P(y)}\}
$$
for $t_1,t_2\in \Phi(P(y))$. For this partition, we have the following lemma.
\begin{lemma}
 $ S_y(a,b,t_1,t_2)$ is either empty or an arithmetic progression and we have
 $$\mathbf{d}S_y(a,b,t_1,t_2)=\begin{cases}0&\ \it{if} P(y)\nmid r_1-r_2-at_1+bt_2,
\\ \frac{k}{abP(y)}&\ \it{if} P(y)\mid r_1-r_2-at_1+bt_2.
\end{cases}$$
\end{lemma}

\proof We may assume $\gcd(a,b)\mid r_1-r_2$, otherwise $S_y(a,b)=\emptyset.$ For $n\in S_y(a,b,t_1,t_2)$, we have $kn+r_1=ax$ and $ kn+r_2=by$ for some $x,y \in \mathbb{Z}$. Therefore, we obtain the linear Diophantine equation
\begin{equation}\label{df}
ax-by =r_1-r_2.
\end{equation}
Writing the congruence conditions as
$$
x=t_1+x'P(y),\ \ \ y=t_2+y'P(y) , \ \ \ x', y' \in \mathbb{Z}.
$$
It follows that
\begin{equation}\label{ap}
aP(y)x'-bP(y)y'=r_1-r_2-at_1+bt_2 .
\end{equation}
If $P(y)\nmid r_1-r_2-at_1+bt_2$, this equation has no solutions, which means that $\mathbf{d}S_y(a,b,t_1,t_2)=0$. When $P(y)\mid r_1-r_2-at_1+bt_2$,  write $P(y)l=r_1-r_2-at_1+bt_2$, which implies that
\begin{equation}\label{x,}
    ax'-by'=l.
\end{equation}
The general solution of \eqref{x,} is $x'=x_0l+mb$ and $y'=y_0l+ma$, $m\in Z$, where $x=x_0$, $y=y_0$ is a particular solution for \eqref{df}. In summary, we have $n\in S_y(a,b,t_1,t_2)$ has the form

$$
kn+r_1=a(t_1+P(y)lx_0)+abP(y)m,
$$
$$
kn+r_2=b(t_2+P(y)ly_0)+abP(y)m,
$$
and any choice of $m$ such that $n\in \mathbb{N}$ puts $n$ in $ S_y(a,b,t_1,t_2).$
Hence, $S_y(a,b,t_1,t_2)$ is an arithmetic progression when nonempty with $$\mathbf{d}S_y(a,b,t_1,t_2)= \frac{k}{abP(y)}.$$
This ends the proof. \qed

We will determine $\mathbf{d}S_y(a,b)$ by counting the number of ordered pairs $(t_1,t_2)$ satisfying $P(y)\mid r_1-r_2-at_1+bt_2.$ For each prime $p \mid P(y)$, we consider the valid pairs $(t_1,t_2)$ modulo $p$ in four cases.

\noindent \textbf{Case 1:} $(p,r_1-r_2)=1$ and $p\nmid ab.$

\noindent In this case, $p\mid r_1-r_2-at_1+bt_2$ if and only if $t_2\equiv b^{-1}(at_1-r_1+r_2)\pmod{p}.$ For each $t_1\in \Phi (P(y))$, we failed to get a valid $t_2\in \Phi (P(y))$ only if $t_1\equiv (r_1-r_2)a^{-1}\pmod p.$
Thus, there are $p-2$ valid ordered pairs modulo $p$.

\noindent \textbf{Case 2:} $(p,r_1-r_2)=1$ and $p\mid ab.$

\noindent
Hence, we cannot have both $p\mid a$ and $p\mid b$, otherwise $p\mid r_1-r_2$, which is a contradiction.
If $p\mid a$, then  $p\mid r_1-r_2-at_1+bt_2$ if and only if $t_2\equiv -(r_1-r_2)b^{-1}\pmod p$,
so $t_1$ is free and $t_2$ is completely determined modulo $p$. Therefore, we have $p-1$
ordered pairs modulo $p$. Similarly, we also have $p-1$ ordered pairs modulo $p$ when $p\mid b.$

\noindent \textbf{Case 3:} $p\mid r_1-r_2$ and $p\nmid ab.$

\noindent
In this case, $p\mid r_1-r_2-at_1+bt_2$ if and only if $t_1\equiv a^{-1}bt_2\pmod {p}.$ Hence,  $t_2$ is free and $t_1$ is completely determined modulo $p$. Then we have $p-1$ ordered pairs modulo $p$.

\noindent \textbf{Case 4:} $p\mid r_1-r_2$ and $p\mid ab.$

\noindent Thus, it is easy to deduce that $p\mid a$ and $p\mid b$ by $p\mid r_1-r_2-at_1+bt_2.$ Hence, $t_1$ and $t_2$ are free and we have $(p-1)^2$ ordered pairs $(t_1,t_2)$ modulo $p.$

Combining the above and using the Chinese reminder theorem, we obtain that
\begin{align*}
\# \{(t_1,t_2)\in \Phi(P(y))^2:P(y)\mid r_1-r_2-at_1+bt_2\}=&\prod_{\substack{p|P(y)\\p\nmid r_1-r_2\\p\nmid ab}}(p-2)\prod_{\substack{p|P(y)\\p\nmid r_1-r_2\\p\mid ab}}(p-1)\\
&\prod_{\substack{p|P(y)\\p\mid r_1-r_2\\p\nmid ab}}(p-1)\prod_{\substack{p|P(y)\\p\mid r_1-r_2\\p\mid ab}}(p-1)^2.
\end{align*}
Hence,
$$
\mathbf{d}S_y(a,b)=\frac k{ab}\prod_{\substack{p|P(y)\\p\nmid r_1-r_2\\p\nmid ab}}\l(1-\frac 2p\r)\prod_{\substack{p|P(y)\\p\nmid r_1-r_2\\p\mid ab}}\l(1-\frac 1p\r)\prod_{\substack{p|P(y)\\p\mid r_1-r_2\\p\nmid ab}}\l(1-\frac 1p\r)\prod_{\substack{p|P(y)\\p\mid r_1-r_2\\p\mid ab}}\frac{ (p-1)^2}p .
$$

The following lemma is from M. Kobayashi and T. Trudgian \cite[Lemma 1]{KT20}.
\begin{lemma}\label{lemma1}
Let
$$g(n):=\l(\frac {\sigma (n)}n \r)^s, \quad
\rho(p^\alpha):=g(p^\alpha)-g(p^{\alpha-1}), $$
and
$$\Lambda_k(s):=\prod_{p\nmid k}\l(1+\frac {\rho(p)} p+\frac {\rho(p^2)} {p^2}+\cdots \r).$$
If $g$ and $h$ are given coprime positive integers with $s\ge 1$ and $x\ge 2,$ then
$$
\sum_{\substack{n\le x\\ n\equiv g\pmod h} }\l(\frac {\sigma (n)}n \r)^s=x\frac {\Lambda_h(s)}h+O((\log h)^s ).
$$
\end{lemma}

Now, we estimate the bounds on $\mathbf{d}B(k,r_1,r_2)$. Since $\mathbf{d}A(k,r_1,r_2)=\mathbf{d}B(k,r_1,r_2)$, we can also get the bounds on $\mathbf{d}A(k,r_1,r_2)$. We say $0$ and $\mathbf{d}S_y(a,b)$ trivial bounds for $\mathbf{d}B_y(a,b)$. For a nontrivial upper bound, we have

\begin{align}
\sum_{\substack{n\in S_y(a,b)\\ n\leq x}}h^s(kn+r_1)&=\sum_{\substack{n\in S_y(a,b)\\  n\in B(k,r_1,r_2)\\n\leq x}}h^s(kn+r_1)+\sum_{\substack{n\in S_y(a,b)\\n\not\in B(k,r_1,r_2)\\ n\leq x}}h^s(kn+r_1)\nonumber\\
&\geq \sum_{\substack{n\in S_y(a,b)\\  n\in B(k,r_1,r_2)\\n\leq x}}h^s(kn+r_2)+\sum_{\substack{n\in S_y(a,b)\\  n\not\in B(k,r_1,r_2)\\n\leq x}}h^s(kn+r_1)\nonumber\\
&\geq h^s(b)|B_y(a,b)\cap [1,x]|+h^s(a)(|S_y(a,b)\cap [1,x]|-|B_y(a,b)\cap [1,x]|).
\end{align}

\
Since $Y(kn+r_2)=b$ , we have $b\mid P(y)$ and $\Lambda_{bP(y)}(s)=\Lambda_{P(y)}(s)$. By Lemma 3.1, the set $S_y(a,b,t_1,t_2)$ is an arithmetic progression when $P(y)\mid r_1-r_2-at_1+bt_2$. By $(t_1,P(y))=1$ and $b\mid P(y)$, we have $(t_1,b)=1.$  Therefore, we obtain that $(t_1+P(y)lx_0,bP(y))=1$. By Lemma 3.3, for such pairs $(t_1,t_2)$ we have
\begin{align*}
\sum_{\substack{n\leq x\\n\in S_y(a,b,t_1,t_2)}}h^s(kn+r_1)&=h^s(a)\sum_{\substack{m\leq (kx+r_1)/a\\ m\equiv t_1+P(y)lx_0 \pmod {bP(y)}}}h^s(m)\\
&=h^s(a)\Lambda_{bP(y)}(s)\frac{kx+r_1}{abP(y)}+O(\log^s x)\\
&=h^s(a)\Lambda_{P(y)}(s)\frac{k}{abP(y)}x+O(\log^s x).
\end{align*}
Summing over all pairs $(t_1,t_2)$, we obtain
\begin{equation}\label{hs}
\sum_{\substack{n\in S_y(a,b)\\n\leq x}}h^s(kn+r_1) \sim h^s(a)\Lambda_{P(y)}(s)\mathbf{d}S_y(a,b)x ,\ \ x\rightarrow \infty.
\end{equation}
Similarly, we have
\begin{equation}\label{hk}
\sum_{\substack{n\in S_y(a,b)\\n\leq x}}h^s(kn+r_2)\sim h^s(b)\Lambda_{P(y)}(s) \mathbf{d}S_y(a,b)x ,\ \ x\rightarrow \infty.
\end{equation}
Hence, dividing both sides of (3.4) by $x$ and taking $x\rightarrow \infty$, by \eqref{hs} we have
$$
h^s(a)\Lambda_{P(y)}(s)\mathbf{d}S_y(a,b)\geq h^s(b)\mathbf{d}B_y(a,b)+h^s(a)\mathbf{d}S_y(a,b)-h^s(a)\mathbf{d}B_y(a,b).
$$
In the case of $h(b)>h(a),$ we arrive at  the upper bound
$$
\mathbf{d}B_y(a,b)\leq \frac {h^s(a)(\Lambda_{P(y)}(s)-1)}{h^s(b)-h^s(a)}\mathbf{d}S_y(a,b).
$$
Observe that $\Lambda_P(s)>1$ for all $s\geq 1,$ this bound is nontrivial when $ h(b)/h(a)>\Lambda_{P(y)}(s)^{1/s},$ which implies that $h(b)>h(a)$.

For a nontrivial lower bound, we have
\begin{align}
\sum_{\substack{n\in S_y(a,b)\\ n\leq x}}h^s(kn+r_2)&=\sum_{\substack{n\in S_y(a,b)\\  n\in B(k,r_1,r_2)\\n\leq x}}h^s(kn+r_2)+\sum_{\substack{n\in S_y(a,b)\\  n\not\in B(k,r_1,r_2)\\n\leq x}}h^s(kn+r_2)\nonumber\\
&\geq \sum_{\substack{n\in S_y(a,b)\\  n\in B(k,r_1,r_2)\\n\leq x}}h^s(kn+r_2)+\sum_{\substack{n\in S_y(a,b)\\  n\not\in B(k,r_1,r_2)\\n\leq x}}h^s(kn+r_1)\nonumber\\
&\geq h^s(b)|B_y(a,b)\cap [1,x]|+h^s(a)(|S_y(a,b)\cap [1,x]|-|B_y(a,b)\cap [1,x]|).
\end{align}
Dividing both sides of (3.7) by $x$ and taking $x\rightarrow \infty$, by \eqref{hk} we have
$$
h^s(b)\Lambda_{P(y)}(s)\mathbf{d}S_y(a,b)\geq h^s(b)\mathbf{d}B_y(a,b)+h^s(a)\mathbf{d}S_y(a,b)-h^s(a)\mathbf{d}B_y(a,b).
$$
If $h(b)<h(a)$, we have
$$
\mathbf{d}B_y(a,b)\geq \frac {h^s(a)-h^s(b)\Lambda_{P(y)}(s)}{h^s(a)-h^s(b)} \mathbf{d}S_y(a,b).
$$
This bound is nontrivial when $h(a)/h(b)>\Lambda_{P(y)}(s)^{1/s},$ which implies $h(b)<h(a)$.

It is easy to verify that
\begin{align*}
\Lambda_{P(y)}(1)&=\prod_{p\nmid P(y)}\l(1+\frac {\rho(p)} p+\frac {\rho(p^2)} {p^2}+\cdots \r)\\
&=\prod_{p\nmid P(y)}\l(1+\frac 1 {p^2}+\frac 1{p^4} +\cdots \r)\\
&=\zeta(2)\prod_{p\mid P(y)}\l(1-\frac 1{p^2}\r).
\end{align*}
Thus, for upper bounds $\Lambda_{P(y)}^+(r)$ for $\Lambda_{P(y)}(r)$, if $s=1$ we use
$$\Lambda^+_{P(y)}(1)=\zeta(2)\prod_{p\mid P(y)}\l(1-\frac 1{p^2}\r).$$
When $s\geq 2$, we use the work of Del\'eglise's  \cite{D97} where we have taken $65536$ to be the maximum prime bound. We obtain
$$\Lambda_{P(y)}^+(s)=\prod_{\substack{p\ prime\\ y<p<65536}}\l(1+\frac {(1+1/p)^s-1} p+\frac s{(p^4-p^2)(1-\frac 1p)^{s-1}} \r)\exp(1.6623114\times 10^{-6}s).
$$
Combining the above, we use the following bounds for $\mathbf{d}B_y(a,b)$:
$$\mathbf{d}B_y(a,b)\geq \mathbf{d}B_y^{-}(a,b)=\begin{cases}\frac{h^s(a)-h^s(b)\Lambda_{P(y)}^+(s)}{h^s(a)-h^s(b)}\mathbf{d}S_y(a,b)&\t{if}\ h(a)/h(b)>\Lambda_{P(y)}^+(s)^{1/s},
\\0&\t{otherwise},\end{cases}$$
$$\mathbf{d}B_y(a,b)\leq\mathbf{d}B_y^{+}(a,b)=\begin{cases}\frac{h^s(a)(\Lambda_{P(y)}^+(s)-1)}{h^s(b)-h^s(a)}\mathbf{d}S_y(a,b)&\t{if}\ h(b)/h(a)>\Lambda_{P(y)}^+(s)^{1/s},
\\\mathbf{d}S_y(a,b)&\t{otherwise}.\end{cases}$$
And then we  give a estimate on the natural density of $A(k,r_1,r_2)$:
$$
\sum_{a,b\in S(y)} \mathbf{d}B_y^{-}(a,b)\leq\mathbf{d}B(k,r_1,r_2)= \mathbf{d}A(k,r_1,r_2)\leq \sum_{a,b\in S(y)} \mathbf{}B_y^{+}(a,b).
$$
\section{Proofs of Theorem 1.2 and Theorem 1.3 }
\setcounter{lemma}{0}
\setcounter{theorem}{0}
\setcounter{equation}{0}
\setcounter{conjecture}{0}
\setcounter{remark}{0}
\setcounter{corollary}{0}

\proof We complete the proofs of Theorem 1.2 and Theorem 1.3 by computer.

 For $k=3$, $r_1=2$ and $r_2=0,$ we have
$S_y(a,b)=\emptyset $ unless $3\nmid a$, $3\mid b$ and $\gcd(a,b)\mid 2$. We fix the parameters $y,z$ and $s_{max}$. Then we recursively run through $a\in S(y)\cap [1,z]$ with $3\nmid a$. For each $a$ we also recursively run through  $b\in S(y)\cap [1,z/a]$ with $3\mid b$. For a given pair $(a,b)$, we calculate $\mathbf{d}B_y^{+}(a,b)$ and $\mathbf{d}B_y^{-}(a,b)$ for $1\leq s\leq s_{max}$. Thus, by the estimate given in Section 3, we obtain the bounds on $\mathbf{d}A(3,2,0) $. Changing $y,z$ and $s_{max}$, we may get more precise bounds on $\mathbf{d}A(3,2,0) $. We take $s_{max}=3$, $y=67$, $z=3*10^7$  for the calculation of lower bound and take $s_{max}=3$, $y=89$, $z=9*10^{10}$ for the calculation of upper bound. Then we get
$$ 0.0361784 \leq \mathbf{d}A(3,2,0) \leq 0.1784994.  $$
For $k=4$, $r_1=1$ and $r_2=0,$ we have
$S_y(a,b)=\emptyset $ unless $a$ is odd, $4\mid b$ and $\gcd(a,b)=1$. Similarly, we calculate the bounds on $\mathbf{d}A(4,1,0)$. We take $s_{max}=3$, $y=37$, $z=9*10^8$ for the calculation of lower bound and take $s_{max}=8$, $y=41$, $z=2*10^{10}$ for the calculation of upper bound. Thus, we obtain
$$ 0.0107553 \leq \mathbf{d}A(4,1,0) \leq 0.0207187.  $$

All the calculations were done by Mathematica. We notice that the accuracy of our bounds is affected by the estimate of $\Lambda_{P(y)}(s)$ and the restriction of $a,b\in S(y)$.

The above two results correspond to the entries A398676\cite{oeis76} and A398677\cite{oeis77} in the OEIS, edited by A.Eldar.


\begin{thebibliography}{99}

\bibitem{B51} N. G. de Bruijn, \textit{On the number of positive integers $\leq$ x and free of prime factors >y.}, \emph{Nederl. Akad. Wetensch. Proc. Ser. A}, {\bf 54} (1951), 50--60.
\bibitem{D97}  M. Del\'eglise, \textit{Bounds for the density of abundant integers}, \emph{Exp.Math}, {\bf 7(2)} (1997), 137--143.
\bibitem{oeis76}
A. Eldar, ``Sequence A398676,'' in \textit{The On-Line Encyclopedia of Integer Sequences}, OEIS Foundation Inc., 2026, \url{https://oeis.org/A398676}.
\bibitem{oeis77}
A. Eldar, ``Sequence A398677,'' in \textit{The On-Line Encyclopedia of Integer Sequences}, OEIS Foundation Inc., 2026, \url{https://oeis.org/A398677}.
\bibitem{E34}  P. Erd\"os, \textit{On the Density of the Abundant Numbers}, \emph{J. London Math. Soc.} {\bf 9} (1934), 278-–282.
\bibitem{E36}  P. Erd\"os, \textit{On a problem of Chowla and some related problems}, \emph{Proc. Camb. Philos. Soc.} {\bf 32} (1936), 530--540.
\bibitem{E57}  P. Erd\"os, \textit{Some unsolved problems}, \emph{Michigan Math. J.} {\bf 4} (1957), 291--300.
\bibitem{E80}  P. Erd\"os, K. Gy\"ory, and Z. Papp, \textit{On some new properties of functions $\sigma(n),\phi(n)$ ,d(n) and v(n)} \emph{Mat.Lapok}, {\bf 28} (1980), 125--131.
\bibitem{EP87}  P. Erd\"os, C. Pomercance and A.S\'ark\"ozy, \textit{On locally Repeated values of certain arithmetic functions. II}, \emph{Acta Math. Hungarica}, {\bf 49} (1987), 251--259.
\bibitem{G13}  T. H. Gr\"onwall, \textit{Some asymptotic expressions in the theory of numbers}, \emph{ Trans. Amer. Math. Soc.}, {\bf 14} (1913), 113--122.
\bibitem{H86} A. Hildebrand, {\it On the number of positive integers $\geq$ x and free of prime factors >y}, \emph{ J. Number Theory}, {\bf 22} (1986), 289--307.
\bibitem{KT20}  M. Kobayashi and T. Trudgian, \textit{On integers of which $\sigma(2n+1)\geq \sigma(2n)$}, \emph{ J. Number Theory}, {\bf 215} (2020), 138--148.
\bibitem{K14}  M. Kobayashi, \textit{A new series for the density of abundant numbers}, \emph{ Int. J. Number Theory}, {\bf 10(1)} (2014), 73--84.
\bibitem{L87}  M. Laub, \textit{Problems and Solutions: Advance Problem: 6555}, \emph{Acta Math. Monthly}, {\bf 94(8)} (1987), 800.
\bibitem{LP15}  F. Luca and C. Pomercance, \textit{The range ot the sum-of-proper-divisors functions}, \emph{Acta Arith.}, {\bf 168(2)} (2015), 187--199.
\bibitem{M90}  L. E. Mattics., \textit{Problems and Solutions:Solutions of Advanced Problems}, \emph{Amer.Math.Monthly}, {\bf 97(4)} (1990), 351–-353.
\bibitem{P77}  C. Pomercance, \textit{On the distribution of amicable numbers. II}, \emph{Journal für die reine und angewandte Mathematik (Crelles Journal)}, {\bf 325} (1977), 183--388.
\bibitem{S73}  H. N. Shapiro, \textit{Addition of functions in probabilistic number theory}, \emph{Comm. Pure Appl. Math}, {\bf 26(1)} (1973), 55-84.
\bibitem{WC22}  R.-J. Wang and Y.-G. Chen, \textit{On positive integers m with $\sigma_l(2n+1)<\sigma_l(2n)$}, \emph{Periodica Mathematica Hungarica}, {\bf 85} (2022), 210--224.

\end{thebibliography}
\end{document}